\documentclass[microtype]{jloganal}
\usepackage{amsfonts}
\usepackage[numbers]{natbib}
\usepackage[fit]{truncate}
\usepackage{amsthm}

\usepackage[UKenglish]{babel}

\usepackage{verbatim}
\usepackage[utf8x]{inputenc}

\usepackage{pinlabel}

\title[Value Groups and Residue Fields of Models of Real Exponentiation]{Value Groups and Residue Fields of\\  Models of Real Exponentiation}
\author[L\,S Krapp]{Lothar Sebastian Krapp}
\givenname{Lothar Sebastian}
\surname{Krapp}

\address{Fachbereich Mathematik und Statistik\\
	Universität Konstanz
	\\ \newline 78457 Konstanz
	\\Germany}
\email{sebastian.krapp@uni-konstanz.de}
\urladdr{http://www.math.uni-konstanz.de/~krapp/}

\arxivreference{1803.03153}
\arxivpassword{4irvn}

\volumenumber{}
\issuenumber{}
\publicationyear{}
\papernumber{}
\startpage{}
\endpage{}
\doi{}
\MR{}
\Zbl{}
\received{}
\revised{}
\accepted{}
\published{}
\publishedonline{}
\proposed{}
\seconded{}
\corresponding{}
\editor{}
\version{}
\numberwithin{equation}{section}

\theoremstyle{plain}
\newtheorem{theorem}{Theorem}[section]
\newtheorem{corollary}{Corollary}[section]
\newtheorem{lemma}{Lemma}[section]
\newtheorem{proposition}{Proposition}[section]
\newtheorem{question}{Question}[section]
\newtheorem{conjecture}{Conjecture}[section]
\theoremstyle{definition}

\newtheorem{construction}{Construction}[section]
\newtheorem{remark}{Remark}[section]

\makeatletter
\let\c@corollary=\c@theorem
\let\c@lemma=\c@theorem
\let\c@proposition=\c@theorem
\let\c@question=\c@theorem
\let\c@conjecture=\c@theorem
\let\c@definition=\c@theorem
\let\c@construction=\c@theorem
\let\c@remark=\c@theorem
\makeatother

\makeautorefname{theorem}{Theorem}
\makeautorefname{corollary}{Corollary}
\makeautorefname{lemma}{Lemma}
\makeautorefname{proposition}{Proposition}
\makeautorefname{question}{Question}
\makeautorefname{conjecture}{Conjecture}
\makeautorefname{defintion}{Definition}
\makeautorefname{construction}{Construction}
\makeautorefname{remark}{Remark}

\newcommand{\N}[0]{\mathbb{N}}
\newcommand{\e}[0]{\mathrm{e}}
\newcommand{\EXP}[0]{\mathrm{EXP}}
\newcommand{\expres}[0]{\overline{\exp}}
\newcommand{\U}{\mathcal{U}}

\newcommand{\setbr}[1]{\left\{ #1 \right\}}
\newcommand{\pow}[1]{\!\left(\!\left( #1 \right)\!\right)}

\renewcommand{\L}[0]{\mathcal{L}}
\newcommand{\Lexp}[0]{\mathcal{L}_{\mathrm{exp}}}
\newcommand{\Fexp}[0]{\mathcal{F}_{\mathrm{exp}}}
\newcommand{\Rexp}[0]{\mathbb{R}_{\mathrm{exp}}}

\newcommand{\Trcf}[0]{T_{\mathrm{RCF}}}
\newcommand{\Kres}[0]{\overline{K}}
\newcommand{\Kexp}[0]{\mathcal{K}_{\mathrm{exp}}}
\newcommand{\Kexpres}[0]{\overline{\mathcal{K}}_{\overline{\mathrm{exp}}}}
\newcommand{\Texp}[0]{T_{\mathrm{exp}}}
\newcommand{\expl}[0]{\exp_{\mathrm{L}}}
\newcommand{\expm}[0]{\exp_{\mathrm{M}}}
\newcommand{\expr}[0]{\exp_{\mathrm{R}}}
\newcommand{\ol}[1]{\overline{#1}}
\newcommand{\ul}[1]{\underline{#1}}
\newcommand\restr[2]{{
		\left.\kern-\nulldelimiterspace 
		#1
		\vphantom{\big|} 
		\right|_{#2}
}}

\DeclareMathOperator{\dcl}{dcl}
\DeclareMathOperator{\supp}{supp}

\def\co{\colon\thinspace}

\begin{document}	
	
	\begin{abstract}
		Let $F$ be an archimedean field, $G$ a divisible ordered abelian group and $h$ a group exponential on $G$. A triple $(F,G,h)$ is realised in a non-archimedean exponential field $(K,\exp)$ if the residue field of $K$ under the natural valuation is $F$ and the induced exponential group of $(K,\exp)$ is $(G,h)$. We give a full characterisation of all triples $(F,G,h)$ which can be realised in a model of real exponentiation in the following two cases: i) $G$ is countable. ii) $G$ is of cardinality $\kappa$ and $\kappa$-saturated for an uncountable regular cardinal $\kappa$ with $\kappa^{<\kappa} = \kappa$. Moreover, we show that for any o-minimal exponential field $(K,\exp)$ satisfying the differential equation $\exp' = \exp$, its residue exponential field is a model of real exponentiation.
	\end{abstract}

	\begin{asciiabstract}
	Let F be an archimedean field, G a divisible ordered abelian group and h a group exponential on G. A triple (F,G,h) is realised in a non-archimedean exponential field (K,exp) if the residue field of K under the natural valuation is F and the induced exponential group of (K,exp) is (G,h). We give a full characterisation of all triples (F,G,h) which can be realised in a model of real exponentiation in the following two cases: i) G is countable. ii) G is kappa-saturated and of cardinality kappa for an uncountable regular cardinal kappa with kappa^( < kappa) = kappa. Moreover, we show that for any o-minimal exponential field (K,exp) satisfying the differential equation exp' = exp, its residue exponential field is a model of real exponentiation.
    \end{asciiabstract}

	\maketitle


	\section{Introduction}
	The following is well-known: Given an archimedean real closed field $F$ and a divisible ordered abelian group $G$, there exists a non-archimedean real closed field $K$ such that the residue field of $K$ under the natural valuation $v$ is $F$ and its value group is $G$. For instance, one can take $K$ to be the ordered field of generalised power series $F\pow{G}$. Due to Tarski's quantifier elimination result in the ordered field of real numbers $(\R,+,\cdot,0,1,<)$ (cf. \cite{tarski}) and the resulting completeness of $\Trcf$, the theory of real closed fields, one can rephrase the above fact as follows: For any archimedean real closed field $F$ and any divisible ordered abelian group $G$ there exists a non-archimedean model $K$ of $\Trcf$ such that the residue field of $K$ is $F$ and its value group under $v$ is $G$.
	
	In the spirit of Tarski's work, we can ask a similar question once we expand the ordered field of real numbers by its standard exponential function $(\R,+,\cdot,0,1,<,\exp)$ and consider the theory of this structure $\Texp$. 
	In general, an exponential function on a non-archimedean field induces a richer structure on its value group and results in an interplay between its residue field, its value group and the rank of its value groups. This interplay is vastly explored in S~Kuhlmann \cite{kuhlmann}. Our main question is the following.
	
	\begin{question}\label{qu3} 
		Let $F$ be an archimedean real closed field, $G$ a divisible ordered abelian group and $h$ a group exponential on $G$. What are necessary and sufficient conditions on $F$, $G$ and $h$ such that there exists a non-archimedean model $(K,+,\cdot,0,1,<,\exp)$ of $\Texp$ such that its residue field under $v$ is $F$ and its induced exponential group is $(G,h)$? 
	\end{question}
	
	S~Kuhlmann posed two further questions which are strongly related to the context of \fullref{qu3} and which also inspire parts of this work.
	
	\begin{question}\label{qu1} 
		Under what conditions is an exponential group the value group of an ordered exponential field? 
	\end{question}
	
	\begin{question}
		\label{qu2} Does the value group of any ordered exponential field have the lifting property? 
	\end{question}

	Both questions are answered in \cite[Chapter~1]{kuhlmann} in the case that the group under consideration is countable. Based on this and some model theoretic properties of $\Texp$ such as model completeness and o-minimality (cf. Wilkie \cite{wilkie}), we will give a complete characterisation of triples $(F,G,h)$ which can be realised in a non-archimedean model of $\Texp$ under the condition that $G$ is countable.
	
	For a given regular uncountable cardinal $\kappa$,  Kuhlmann--Shelah \cite{kuhlmann2} provides a construction method of $\kappa$-bounded Hahn sums which are the value groups of models of $\Texp$. Based on this and further assuming that $\kappa^{<\kappa} = \kappa$, we will answer all three questions for $\kappa$-saturated divisible ordered abelian groups of cardinality $\kappa$.
	
	The paper is structured as follows. In \fullref{secprelim} we will set up the terminology, notation and some facts we need throughout the other sections. In \fullref{secomin} we prove some general results on o-minimal exponential fields in connection to their residue fields and exhibit the links to the open decidability problem of the real exponential field. Moreover, we give a characterisation of all archimedean fields which can be realised as the residue field of models of real exponentiation. 
	In \fullref{secctbl} we provide answers to the questions above in the case that the group under consideration is countable. This raises the question whether there exist uncountable models of real exponentiation with countable value group. A construction method of such models is given at the end of the section. We continue by considering $\kappa$-bounded power series fields for a regular uncountable cardinal $\kappa$ with $\kappa^{<\kappa}=\kappa$ in \fullref{seckappasat}. This leads to answers to the questions above for $\kappa$-saturated groups of cardinality $\kappa$. Finally, in \fullref{seccontr}, we point out how the results of \fullref{secctbl} and \fullref{seckappasat} can be used to give a characterisation of triples $(F,G,\chi)$ which can be realised in a model of real exponentiation, where $\chi$ is a contraction on $G$ rather than a group exponential.

	\section{Preliminaries}\label{secprelim}
	
	In this section we will outline the algebraic and model theoretic preliminaries on ordered exponential fields and exponential groups. We follow the notation from S~Kuhlmann \cite{kuhlmann}, which also contains further details.
	
	Let $(K,+,\cdot,0,1,<)$ be an ordered field. An {\bf exponential} $\exp$ on $K$ is an order-preserving isomorphism from the ordered additive group $(K,+,0,<)$ to the ordered multiplicative group $(K^{>0},\cdot,1,<)$. The structure $(K,+,\cdot,0,1,\allowbreak <,\exp)$ is then called an {\bf ordered exponential field} and denoted by $\Kexp$. The language of ordered exponential fields is denoted by $\Lexp = (+,\cdot,0,1,\allowbreak <,\exp)$.
	If it is clear from the context, we only write $K$ for the ordered field $(K,+,\cdot,0,1,<)$ and also $(K,\exp)$ for the ordered exponential field $(K,+,\cdot,0,1,<,\exp)$. 
	
	The most prominent example of an ordered exponential field is $\Rexp = (\R,\exp_\R)$, where $\exp_\R$ is the standard exponential function $\exp_\R\co x \mapsto \e^x$ on $\R$. We denote the complete theory of $\Rexp$ by $\Texp$. If it is clear from the context, we also denote the standard exponential on $\R$ by $\exp$. If a real closed field $K$ admits an exponential $\exp$ such that $(K,\exp) \models \Texp$, then $\exp$ is called a \textbf{$\Texp$-exponential} on $K$ and $(K,\exp)$ a \textbf{model of real exponentiation}.
	
	We denote by $v$ the \textbf{natural valuation} on $K$, by $R_v$ the \textbf{valuation ring}, by ${I}_v$ the \textbf{valuation ideal}, by $\mathcal{U}_v^{>0} = \{x \in K^{>0} \mid v(x) = 0\}$ the \textbf{multiplicative group of positive units}, and by $\overline{K} = R_v / {I}_v$ the \textbf{residue field} of $K$. For $a \in R_v$, we denote by $\ol{a}$ its residue in $\Kres$. If not stated otherwise, we denote the value group $v(K)$ of $K$ by $G$, the natural valuation on $G$ by $v_G$ and the value set $v_G(G)$ of $G$ by $\Gamma$. 
	
	For any ordered group $H$ and any $\gamma \in v_H(H)$, let $H^{\gamma} = \{x \in H  \mid v_H(x)  \geq \gamma\}$ and $H_{\gamma} = \setbr{x \in H \mid v_H(x) > \gamma}$. Then $B(H,\gamma) = H^\gamma / H_\gamma$ is the \textbf{archimedean component} corresponding to $\gamma$. 
	
	An exponential $\exp$ on $K$ is called \textbf{$v$--compatible} if $\exp(R_v) = \mathcal{U}_v^{>0}$ and $\exp(I_v) = 1 + I_v$. 
	A $v$--compatible exponential on $K$ induces a \textbf{residue exponential} $\expres$ on the residue field $\Kres$ given by $\expres(\ol{a}) = \ol{\exp(a)}$. We denote $(\Kres,\expres)$ by $\Kexpres$ and call it the \textbf{residue exponential field} of $\Kexp$.
	
	For any ordered field $K$ there exists an \textbf{additive lexicographic de\-com\-po\-si\-tion} of the ordered additive group $(K,+,0,<)$, which we denote by $$K = \mathbf{A} \amalg \mathbf{A}' \amalg I_v,$$ where $\mathbf{A}$ and $\mathbf{A}'$ are unique up to isomorphism, $\mathbf{A}'$ is order isomorphic to $(\Kres,+,0,<)$, $v(\mathbf{A}) = G^{<0}$ and the non-zero archimedean components of $\mathbf{A}$ and ${I_v}$ are all isomorphic to $(\Kres,+,0,<)$. If, moreover, the ordered multiplicative group of positive elements $(K^{>0},\cdot,1,<)$ is divisible, then there exists a \textbf{multiplicative lexicographic decomposition} of $(K^{>0},\cdot,1,<)$, denoted by $$K^{>0} = \mathbf{B} \amalg \mathbf{B}' \amalg (1 + I_v),$$ where $\mathbf{B}$ is order isomorphic to $G$ via $-v$ and $\mathbf{B}'$ is order isomorphic to $(\Kres^{>0},\cdot,1,<)$.
	
	Let $\exp$ be a $v$--compatible exponential on $K$. For a given additive lexicographic decomposition $K = \mathbf{A} \amalg \mathbf{A}' \amalg I_v$, we can set $\mathbf{B} = \exp(\mathbf{A})$ and $\mathbf{B}' = \exp(\mathbf{A}')$. Then $\mathbf{B} \amalg \mathbf{B}' \amalg (1 + I_v)$ is a multiplicative lexicographic decomposition of $K^{>0}$. We denote by $\expl$ the restriction of $\exp$ to $\mathbf{A}$ and call it \textbf{$v$--left exponential}. The restriction of $\exp$ to $\mathbf{A}'$ is denoted by $\expm$ and the restriction to $I_v$ by $\expr$. These are called\textbf{ $v$--middle} and \textbf{$v$--right exponentials}, re\-spec\-tive\-ly. 
	If $K^{>0}$ is a divisible multiplicative group and both an additive and a multiplicative lexicographic decomposition are given, then any isomorphisms from $\mathbf{A}$ to $\mathbf{B}$, from $\mathbf{A}'$ to $\mathbf{B}'$ and from $I_v$ to $1+I_v$ are also called $v$--left, $v$--middle and $v$--right exponentials on $K$, respectively. For given $v$--left, $v$--middle and $v$--right exponentials $\expl$, $\expm$ and $\expr$ on $K$, we can glue them together to obtain an exponential on $K$ given by 
	\begin{align*}
		\exp = \expl \amalg \expm \amalg \expr\co K \to K^{>0}, \\
		a + a' + \varepsilon \mapsto
		\expl(a)\expm(a')\expr(\varepsilon),
	\end{align*} 
	where $a \in \mathbf{A}$, $a' \in \mathbf{A}'$ and $\varepsilon \in I_v$.
	
	Let $e_\R = \restr{\exp_\R}{[0,1]}$, the {\bf restricted exponential function} on $\R$, and denote the complete theory of $(\R,e_\R)$ by $T_e$. If it is clear from the context, we also write $e$ instead of $e_\R$.
	In the language $\Lexp$ we can state the {\bf growth axiom scheme (GA)}, which consists of the collection of the following sentences (one for each natural number $n\geq 1$):
	\[\forall x \ (x \geq n^2 \to \exp(x) > x^n).\]
	A $v$--left exponential $\expl$ on an ordered field $K$ is called a \textbf{(GA)--$v$--left exponential}, if it satisfies (GA) on the compatible complement $\mathbf{A}$ to $R_v$, i.e. for any $x \in \mathbf{A}$ we have $\exp(x) > x^n$. Due to Ressayre \cite{ressayre}, any ordered field $K$ which admits both a unary function $e$ with $(K,e) \models T_e$ and a (GA)--$v$--left exponential $\expl$, also admits a $\Texp$--exponential $\exp$ extending both $e$ and $\expl$.
	
	Now let $G$ be an ordered abelian group. An isomorphism of chains $h: \Gamma \to G^{<0}$ is called a \textbf{group exponential} on $G$ and the pair $(G,h)$ an \textbf{exponential group}. 	
	Suppose that there is an archimedean abelian group $A$ such that for any $\gamma \in \Gamma$ we have $B(G,\gamma) \cong A$. Then $(G,h)$ is called an \textbf{exponential group in $A$}. If for any $g \in G^{<0}$ we have $h(v_G(g)) > g$, then $(G,h)$ is called a \textbf{strong exponential group}. 
	
	Consider a $v$--compatible exponential field with additive and mul\-ti\-pli\-ca\-tive lexicographic decompositions as above. Then $-v\circ\expl$ induces an isomorphism from $\mathbf{A}$ to $G$, which we denote by $\tilde{\ell}_{\exp}$. Set $\tilde{h}_{\exp} = \tilde{\ell}_{\exp}^{-1}$. Taking valuations on either side, we obtain an isomorphism of chains $h_{\exp}$ given by
	$$h_{\exp}\co \Gamma \to G^{<0}, v_G(g) \mapsto v(\tilde{h}_{\exp}(g)).$$ The exponential group $(G,h_{\exp})$ in $(\Kres,+,\ol{0},<)$ is called the \textbf{exponential group induced by $(K,\exp)$}. The following result connects the notions of a (GA)--$v$--left exponential and a strong exponential group.
			
	\begin{lemma}{\rm \cite[Corollary~2.6]{kuhlmann}}\label{lemkuhl20} \qua
		Let $\Kexp$ be an ordered exponential field, $\expl$ the induced $v$--left exponential and $(G,h_{\exp})$ the induced exponential group. Then $\expl$ is a (GA)--$v$--left exponential if and only if the induced exponential group $(G,h_{\exp})$ is a strong exponential group.
	\end{lemma}
	
	For an archimedean field $F$ and an exponential group $(G,h)$, we call the triple $(F,G,h)$ an \textbf{exponential valuation triple}. Such a triple is realised in a $v$--compatible  non-archimedean field $\Kexp$ if $\Kres = F$ and $\Kexp$ induces the exponential group $(G,h)$.
	
	Let $H$ be an ordered abelian group. An automorphism $\tau$ on $H$ induces an au\-to\-mor\-phism $\tilde{\tau}$ on $v_H(H)$ by $$\tilde{\tau}(v_H(h)) = v_H(\tau(h)).$$ We say that $H$ has the \textbf{lifting property} if any automorphism $\sigma$ on $v_H(H)$ lifts to an automorphism $\tau$ on $H$, i.e. $\tilde{\tau} = \sigma$.
	
	Let $\Gamma$ be an ordered set and $H$ an archimedean abelian group. The \textbf{Hahn product} with index set $\Gamma$ and components in $H$ is denoted by $H^\Gamma$. It consists of all maps $g: \Gamma \to H$ with well-ordered \textbf{support} $\supp g = \{\gamma \in \Gamma \mid g(\gamma)\neq 0 \}$. Endowed with lexicographic order and pointwise addition, $H^\Gamma$ is called a \textbf{Hahn group}. 
	For $G = H^\Gamma$ the natural valuation $v_G$ is given by $v_G(g) = \min \supp g$.
	The subgroup of $H^\Gamma$ of all elements with finite support is denoted by $\coprod_\Gamma H$. It is called the \textbf{Hahn sum} with index set $\Gamma$ and components in $H$. For an ordered field $F$, we write $\coprod_\Gamma F$ for the Hahn sum over its ordered additive group $(F,+,0,<)$.
	
	Let $F$ be an archimedean field and $H\neq 0$ an ordered abelian group. We denote by $F\pow{H}$ the \textbf{field of generalised power series} with coefficients in $F$ and exponents in $G$. As ordered group, this is just $F^H$, where $F$ is the ordered abelian group $(F,+,0,<)$ and $H$ the ordered set $(H,<)$. We denote the \textbf{monomials} in $F\pow{H}$ by $t^h$ for any $h\in H$, that is, $t^h$ is the map sending $h$ to $1\in F$ and any other element of $H$ to $0$. Any element in $s\in F\pow{H}$ can be expressed as $s = \sum_{h\in H} s_ht^h$, where $s_h = s(h)$. With this representation we obtain multiplication by $t^h\cdot t^{h'} = t^{h+h'}$ and distributivity. This makes $F\pow{H}$ an ordered field. The natural valuation on $F\pow{H}$ is given by $v(s) = \min \supp s$.
	
	Let $\kappa$ be a regular uncountable cardinal. We denote the \textbf{$\kappa$--bounded Hahn group} by $(H^\Gamma)_\kappa$. This is the subgroup of $H^\Gamma$ consisting of its elements whose support has cardinality less than $\kappa$. Similarly, $F\pow{H}_\kappa$ is the \textbf{$\kappa$--bounded power series field} consisting of the elements $s\in F\pow{H}$ with $|\supp s| < \kappa$.

	\section{O-minimal Exponential Fields}\label{secomin}
	
	A totally ordered structure $(M,<,\ldots)$ is called \textbf{o-minimal} if any parametrically definable subset of $M$ is a finite union of open intervals and points in $M$.
	O-minimal structures expanding a group exhibit a tame geometric behaviour of definable functions and have certain useful model theoretic properties. For instance, Taylor expansions exist for definable functions and every such structure has a unique prime model. These and further geometric and model theoretic properties of definable functions in o-minimal structures can be found in van den Dries \cite{dries} and Pillay--Steinhorn \cite{pillay,pillay2}.
	
	In order to start our investigation of countable exponential groups, we need some model theoretic results on o-minimal exponential fields connected to their residue fields. These structures are also of great interest on their own, as the study of these has connections to open problems such as the decidability of $\Rexp$ and consequently Schanuel's Conjecture. An overview of these connections is given in Krapp \cite{krapp}. We will also point out how the results in this section relate to them. The crucial point in our treatment of o-minimal exponential fields is the result that $\Rexp$ is model complete and o-minimal, which was first proved in Wilkie \cite{wilkie}.
	
	We denote by $\EXP$ the first-order $\Lexp$--sentence stating the differential equation $\exp' = \exp$ with initial condition $\exp(0) = 1$. We call ordered exponential fields satisfying this sentence {\bf $\EXP$--fields} and exponentials satisfying this differential equation {\bf $\EXP$--ex\-po\-nen\-tials}. 
	
	Our first main result will show that the residue exponential field of an o-minimal $\EXP$--field is an elementary substructure of $\Rexp$.
	In order to show $v$--compatibility, we use the following lemma (cf. \cite[Lemma~1.17]{kuhlmann}).
	
	\begin{lemma}\label{lemmavcomp}
		Let $\Kexp$ be an ordered exponential field. Then $\exp$ is $v$--compatible if and only if $v(\exp(1) - 1) = 0$.
	\end{lemma}
	
	\begin{proposition}\label{mainthmscexp}
		Let $\Kexp$ be a definably complete $\EXP$--field. Then $\exp$ is $v$--compatible and $\Kexpres \subseteq \Rexp$.
	\end{proposition}
	
	\begin{proof}
		By the Taylor expansion of $\exp$ in $-1$, we obtain for any odd $n \in \N$ that
		\[\sum_{k=0}^n \frac{(-1)^k}{k!} < \exp(-1) < \sum_{k=0}^{n+1} \frac{(-1)^k}{k!}.\]
		Note that by \fullref{lemmavcomp}, these bounds for suitably large $n$ imply $v$--compatibility, as $\exp(-1) = \exp(1)^{-1}$. Moreover,
		taking residues, this gives us $\expres(-1) = \e^{-1}$ and thus $\expres(1) = \exp_\R(1)$. Now for any $q \in \Q$ we have $\expres(q) = \exp_\R(1)^q = \exp_\R(q)$. Since $\expres$ and $\exp_\R$ agree on $\Q$ and $\expres$ is continuous on $\Kres$ (as the preimage of an open interval under $\expres$ is an open interval or the empty set), we obtain that $\expres$ and $\exp_\R$ must already agree on all of $\Kres$, as required.
	\end{proof} 
	
	For a structure $\mathcal{M}$ and a subset $A \subseteq M$, we denote by $\dcl(A;\mathcal{M})$ the \textbf{definable closure} of $A$ in $\mathcal{M}$. If $\mathcal{M}$ is an o-minimal expansion of an ordered group, we also say that $\dcl(A;\mathcal{M})$ is an elementary substructure of $\mathcal{M}$, meaning that $\dcl(A;\mathcal{M})$ is the domain of the elementary substructure of $\mathcal{M}$ obtained by restricting functions and relations to $\dcl(A;\mathcal{M})$. Moreover, we let $\Z[x_1,\ldots,x_n]_{\exp}$ be the ring in $n$ variables of exponential polynomials over $\Z$ of exponentiation depth $1$, that is, $\Z[x_1,\ldots,x_n]_{\exp}$ consists of expressions of the form $p(x_1,\ldots,x_n,\exp(x_1),\ldots,\exp(x_n))$ for some $p(x_1,\ldots,x_n,y_1,\ldots,y_n) \in \Z[x_1,\allowbreak\ldots, x_n,y_1,\ldots,y_n]$. If the length of a tuple $(x_1,\ldots,x_n)$ is clear from the context, then we abbreviate it by $\ul{x}$.
	
	In the proof of the following theorem, we only require that $\Kexp$ is a definably complete rather than an o-minimal $\EXP$--field. However, by Fornasiero--Servi \cite[Corollary~8.2 \& Remark~8.3]{fornasiero2} any definably complete $\EXP$--field is already o-minimal (where the assumption \emph{definably Baire} is removed in Hieronymi \cite{hieronymi}).\footnote{We thank the anonymous referee for pointing out this fact and giving the essential references.}
	
	\begin{theorem}\label{elsubstrtexp}
		Let $\Kexp$ be an o-minimal $\EXP$--field. Then $\Kexpres \preceq \Rexp$.
	\end{theorem}
	
	\begin{proof}
		Since $\Kexpres \subseteq \Rexp$, by o-minimality of $\Rexp$ it remains to show that $\dcl(\Kres;\allowbreak\Rexp) = \Kres$.
		
		Let $g$ be a definable function on $\Rexp$. By model completeness of $\Rexp$, we can express $g(\ul{x}) = y$ by a formula of the form $$\exists \ul{z} \ p(\ul{x},y,\ul{z}) = 0,$$
		for some $p(\ul{x},y,\ul{z}) \in \Z[\ul{x},y,\ul{z}]_{\exp}$. Let $\ul{a} \in \Kres$ and $b = g(\ul{a})$. We need to show that $b \in \Kres$. 
		For some $\ul{c} \in \R$ and $N \in \N$ we have 
		\begin{gather*}
		 \Rexp \models p(\ul{a},b,\ul{c}) = 0 \wedge |b| + |\ul{c}| < N.\\
		 \tag*{Hence,} \Rexp \models \exists y \exists \ul{z} \ (p(\ul{a},y,\ul{z}) = 0 \wedge |y| + |\ul{z}| < N).
		\end{gather*}
		If $\Kexpres \models \exists y \exists \ul{z} \ (p(\ul{a},y,\ul{z}) = 0 \wedge |y| + |\ul{z}| < 3N)$, then $p(\ul{a},b',\ul{c}') = 0$ for some $b', \ul{c}' \in \Kres \subseteq \Rexp$. By uniqueness of $b$, we obtain $b = b' \in \Kres$, as required. Assume for a contradiction that
		\begin{align}
			\Kexpres \not\models \exists y \exists \ul{z} \ (p(\ul{a},y,\ul{z}) = 0 \wedge |y| + |\ul{z}| < 3N).\label{fctcontr}
		\end{align}
		Fix $\ul{\alpha} \in R_v$ with $\ul{\ol{\alpha}} = \ul{a}$. Consider $|p(\ul{\alpha},y,\ul{z})|$ on the definable compact set $$D=\setbr{(y,\ul{z}) \in K \mid |y| + |\ul{z}| \leq 2N}.$$ By o-minimality of $\Kexp$, the exponential polynomial $|p(\ul{\alpha},y,\ul{z})|$ attains its minimum $\theta$, say, on $D$ in some point $(\gamma,\ul{\zeta}) \in D$. If $\theta \in I_v$, then $|p(\ul{a},\ol{\gamma},\ul{\ol{\zeta}})| = |p(\ul{\ol{\alpha}},\ol{\gamma},\ul{\ol{\zeta}})| = \ol{\theta} = {0}$. This contradicts \eqref{fctcontr}, as $|\ol{\gamma}|+|\ul{\ol{\zeta}}| \leq 2N < 3N$. Hence, $\theta > \frac 1 k$ for some $k \in \N$. This implies that for all $(y,\ul{z}) \in K$ with $|y| + |\ul{z}| \leq 2N$ we have that $|p(\ul{\alpha},y,\ul{z})| > \frac{1}{k}$. By taking residues, we obtain that for all $(y,\ul{z}) \in \Kres$ with $|y| + |\ul{z}| \leq N$ we have that $|p(\ul{a},y,\ul{z})| \geq \frac{1}{k}$.
		Since $|b| + |\ul{c}| < N$, we can choose a sequence of rational numbers $(p_n)$ and a tuple of sequences of rational numbers $(\ul{q_n})$ with $|p_n|+|\ul{q_n}| < N$ for each $n$, and $ \lim_{n \to \infty} p_n = b$ and $\lim_{n \to \infty} \ul{q_n} = \ul{c}$, where the limits are taken in $\R$. Now $|p(\ul{a},y,\ul{z})|$ is a continuous function. For each $n \in \N$ we have $p(\ul{a},p_n,\ul{q_n}) \geq \frac 1 k$. But by continuity we also have $$\lim_{n \to \infty} p(\ul{a},p_n,\ul{q_n}) = p(\ul{a},b,\ul{c}) = 0,$$ which gives us the required contradiction.
	\end{proof}
	
	\begin{remark} \label{texpre}
		\begin{enumerate}
		\item \fullref{elsubstrtexp} shows in particular that any o-minimal archimedean $\EXP$--field is an elementary substructure of $\Rexp$. This result also follows directly from Laskowski--Steinhorn \cite[Theorem~2.10]{laskowski}, where it is shown that any o-minimal expansion of an archimedean group can be elementarily embedded into a structure with domain $\R$ where both the embedding and the structure are unique.
		
		\item It is shown in Baisalov--Poizat \cite{baisalov} that, given an o-minimal structure $\mathcal{K}$ expanding an ordered field, the residue field $\ol{\mathcal{K}}$ of $\mathcal{K}$ with induced structure is weakly o-minimal. Moreover, if $\mathcal{K}$ is $\omega$--saturated, then $\ol{\mathcal{K}}$ is o-minimal and has domain $\R$. We point out that \fullref{elsubstrtexp} does not need the assumption of $\omega$--saturation.
		
		\item \label{texpre2}	Since the natural valuation ring of a model of $\Texp$ is $\Texp$--convex (cf. van den Dries--Lewenberg \cite[(2.8)]{dries3}), we can deduce that for any $\Kexp$ which is not only o-minimal but also a model of $\Texp$ we already have $\Kexpres \preceq \Kexp$ (cf. \cite[Theorem ~2.12]{dries3}).
		\end{enumerate}
	\end{remark}

	As mentioned at the start of this section, we will now point out how \fullref{elsubstrtexp} relates to the decidability problem of $\Rexp$.
	
	The question whether $\Rexp$ is decidable was posed by Tarski in \cite{tarski}. Macintyre and Wilkie proved in \cite{macintyre} their much-celebrated result that under the assumption of Schanuel's Conjecture, $\Rexp$ is decidable, and, in fact, the decidability of $\Rexp$ is equivalent to a statement they call Weak Schanuel's Conjecture (see also Servi \cite[Section~4.6]{servi}). By axiomatising a recursive subtheory $T_\mathrm{omin}$ of $\Texp$ whose models are o-minimal $\EXP$--fields, Berarducci and Servi showed in \cite{berarducci} that the decidability of $\Rexp$ is also a consequence of the following conjecture.
	
	\begin{conjecture}[Transfer Conjecture for $\EXP$--fields]\label{conj1}
		Let $\Kexp$ be an o-minimal $\EXP$--field. Then $\Kexp \equiv \Rexp$.
	\end{conjecture}
	
	The following proposition settles open questions related to \fullref{conj1} in Krapp \cite[Section~3]{krapp}.
	
	\begin{proposition}\label{prop:settling}
		Let $\Kexp$ be an o-minimal $\EXP$--field. Then the following are equivalent:
		\begin{enumerate}
			\item $\Kexp \equiv \Rexp$.
			\item $\Kexpres \equiv \Kexp$.
			\item $\Kexpres \preceq \Kexp$.
		\end{enumerate}
	\end{proposition}
	
	\begin{proof}
		By \fullref{elsubstrtexp}, (2) implies (1). In \fullref{texpre} (\ref{texpre2}), it is explained how (1) implies (3).
	\end{proof}

	We now turn to models of real exponentiation. \fullref{archfieldchar} will give us a full characterisation of all archimedean fields which are the residue fields of models of $\Texp$.
	
	\begin{lemma}\label{lemexpdclres}
		Let $\Fexp \preceq \Rexp$ and let $\Kexp \models \Texp$ be non-archimedean such that $\Fexp \subseteq \Kexp$. Moreover, let $a \in K\setminus\setbr{0}$ with $v(a) \neq 0$. Then for $F' = \dcl(F(a); \Kexp)$ we have $\ol{F'} = {F}$. 
	\end{lemma}
	
	\begin{proof}
		Assume for a contradiction that there exists $b \in F'$ such that $\ol{b} \notin F$. By the exchange property, $F' = \dcl(F(b);\Kexp)$. By \fullref{texpre} (\ref{texpre2}), there exists an elementary embedding $h\co \Kexpres \to \Kexp$. Let $c = h(\ol{b})$. Now $b$ and $c$ determine the same cut over $F$, that is, $F^{<b} = F^{<c}$. By o-minimality, $b$ and $c$ have the same type over $F$. This implies  $\dcl(F(b);\Kexp) = \dcl(F(c);\Kexp)$. But the latter is isomorphic to a substructure of $\Kexpres$ via $h^{-1}$, whereas $\dcl(F(b);\Kexp) = F'$ is non-archimedean, a contradiction.
	\end{proof}

	\begin{remark}\label{rmkreferee}
	   \fullref{lemexpdclres} holds in greater generality: For instance, the same arguments can be applied to a complete and model complete o-minimal theory $T\supseteq \Trcf$ with an archimedean prime model.\footnote{We thank the anonymous referee for presenting the proof of \fullref{lemexpdclres} and pointing out that it holds in greater generality. The proof of \fullref{lemexpdclres} in a previous version of this paper relied on analytical methods of definable functions in o-minimal structures.}
	\end{remark}
	
	We can now state a full characterisation of all archimedean fields which can be realised as the residue field of a non-archimedean model of $\Texp$.
	
	\begin{theorem}\label{archfieldchar}
		Let $F \subseteq \R$ be an archimedean field. Then the following are equivalent:
		\begin{enumerate}
			\item $F$ is closed under $\exp_\R$ and $(F,\exp_\R) \preceq \Rexp$.
			\item There exists a non-archimedean $\Kexp \models \Texp$ with $\Kres = F$.
		\end{enumerate}
	\end{theorem}
	
	\begin{proof}
		That (2) implies (1) is a consequence of \fullref{elsubstrtexp}.
		
		For the converse, applying \fullref{lemexpdclres} to $(F,\exp_{\R})$ gives us a non-archimedean model of $\Texp$ whose residue field is $F$.
	\end{proof}
	
	\section{Countable Groups}\label{secctbl}
	
	\textbf{From this section onwards, all exponentials we consider are assumed to be $v$--compatible.}
	
	In this section we firstly use the results from \cite[Chapter~1]{kuhlmann} leading up to answering \fullref{qu1} and \fullref{qu2} for the case in which the group under consideration is countable. We will then provide answers to \fullref{qu3} for countable groups.
	Note that for an ordered field $K$, the countability of its value group does in general not imply the countability of $K$. We will show that there exist uncountable models of real exponentiation with countable value group.

	\begin{theorem}{\rm \cite[Theorem\,1.42]{kuhlmann}} \qua \label{kuhl1}
		Let $G$ and $A$ be countable divisible ordered abelian groups such that $A$ is archimedean and $G \neq 0$. Then $G$ is an exponential group in $A$ if and only if $G \cong \coprod_{\Q} A$. 
	\end{theorem}

	\begin{lemma}\label{ctblexpgplemma2}
		Let $\Kexp$ be a non-archimedean exponential field. Then $|\Kres| \leq |G|$. In particular, if $|G| = \aleph_0$, then $|\Kres| = \aleph_0$.
	\end{lemma}
	
	\begin{proof}
		Let $\mathbf{A}$ be the group complement to $R_v$ in the additive lexicographic decomposition of $K$ compatible with $\exp$. We have $\mathbf{A} \cong G$ via $(-v) \circ \exp_{\mathrm{L}}$. Moreover, the archimedean components of $\mathbf{A}$ are isomorphic to $\Kres$. Hence, for some fixed $g \in G^{<0}$ we have $$|\Kres| = |\setbr{a \in \mathbf{A} \mid v(a) \geq g}/\setbr{a \in \mathbf{A} \mid v(a) > g}| \leq |\mathbf{A}| = |G|. \proved$$
	\end{proof}
	
	\begin{proposition}\label{kuhl2cor}\label{lemmaundcondctbl}
		Let $\Kexp$ be a non-archimedean exponential field with countable value group $G$. Then $\Kres$ is countable and $G \cong \coprod_{\Q}\Kres$. 
	\end{proposition}
	
	\begin{proof}
		By \fullref{ctblexpgplemma2}, $\Kres$ is countable. Since $G$ is an exponential group in $\Kres$, we obtain by \fullref{kuhl1} that $G \cong \coprod_{\Q}\Kres$.
	\end{proof}

	We thus obtain a postive answer to \fullref{qu2} in the case that $G$ is countable.

	\begin{proposition}
		Let $\Kexp$ be a non-archimedean exponential field such that $G$ is countable. Then $G$ has the lifting property.
	\end{proposition}

	\begin{proof}
	 Let $\varphi\co G \to \coprod_{\Q}\Kres$ be an isomorphism of ordered groups. This induces an isomorphism of chains $\tilde{\varphi}\co \Gamma \to \Q, v_G(g) \mapsto v(\varphi(g))$. Let $\sigma\co \Gamma \to \Gamma$ be an automorphism. Then $\sigma' = \tilde{\varphi}\circ\sigma\circ\tilde{\varphi}^{-1}$ defines an automorphism on $\Q$. Let $$\tau'\co \coprod_{\Q}\Kres \to \coprod_{\Q}\Kres, s \mapsto s \circ (\sigma')^{-1}.$$ This is an automorphism lifting $\sigma'$. Finally, let $\sigma = \varphi^{-1} \circ \tau' \circ \varphi$. This is the required automorphism which lifts $\tau$.
	\end{proof}
	
	We now turn to \fullref{qu1}. 
	
	\begin{theorem}{\rm \cite[Theorem~1.44]{kuhlmann}} \qua\label{kuhla}
		Let $K$ be a countable non-archimedean ordered field which is root closed for positive elements. Suppose that $e$ is an exponential on $\Kres$. Then the following are equivalent:
		\begin{enumerate}
			\item $K$ admits an exponential $\exp$ such that $\expres = e$.
			\item $v(K) \cong \coprod_\Q \Kres$.
		\end{enumerate}
	\end{theorem}

	By \fullref{lemmaundcondctbl}, we have a necessary condition on $G$ for it to be the value group of an ordered exponential field. The following construction gives us the converse.
	
	\begin{construction}\label{constra}
		Let $(F,e)$ be a countable archimedean exponential field and $G = \coprod_\Q F$. Inspired by \cite[Example 1.45]{kuhlmann}, we will construct a countable non-archimedean exponential field $\Kexp$ with $v(K) = G$. 
		
		Consider $K_0 = F(t^g \mid g\in G) \subseteq F\pow{G}$. Note that $F$ is root closed for positive elements, as it admits an exponential. Hence, $F\pow{G}$ is root closed. Take $K$ to be the root closure for positive elements of $K_0$ in $F\pow{G}$. Then $\Kres = F$ and $v(K) = G$. By \fullref{kuhla}, $K$ admits an exponential lifting $e$.
	\end{construction}
	
	We obtain the following answer to \fullref{qu1} in the countable case.
	
	\begin{proposition}\label{ansqu1}
		Let $G$ be a countable exponential group. Then the following are equivalent:
		\begin{enumerate}
			\item There exists a non-archimedean exponential field $\Kexp$ with $v(K) \cong G$.
			\item There exists a countable archimedean exponential field $(F,e)$ such that $G \cong \coprod_\Q F$.
		\end{enumerate}
	\end{proposition}
	
	\begin{proof}
		If (1) holds, then (2) follows from \fullref{lemmaundcondctbl} setting $F = \Kres$ and $e = \expres$.
		\fullref{constra} shows that (2) implies (1)
	\end{proof}
	
	\begin{remark}
		In \fullref{constra}, we can as well start with a real closed field $F$ and take the real closure of $F(t^g \mid g\in G)$ in $F\pow{G}$ to obtain a countable real closed non-archimedean exponential field $\Kexp$. Hence, for a countable exponential group $G$ there exists a non-archimedean real closed exponential field $\Kexp$ with $v(K) \cong G$ if and only if there exists some countable archimedean real closed exponential field $(F,e)$ such that $G \cong \coprod_\Q F$.
	\end{remark}
	
	We now strengthen this result to answer which countable exponential groups can be realised as the value group of a model of $\Texp$.
	
	\begin{proposition}\label{thmclassexpgps}
		Let $G$ be a countable exponential group. Then the following are equivalent:
		
		\begin{enumerate}
			\item There exists a non-archimedean exponential field $\Kexp \models \Texp$ with $v(K) \cong G$.
			\item There exists some countable $\Fexp \preceq \Rexp$ such that $G \cong \coprod_{\Q} F$.
		\end{enumerate}
	\end{proposition}
	
	\begin{proof}
		Suppose that (1) holds. Set $F = \Kres$. By \fullref{kuhl2cor}, we have $G \cong \coprod_{\Q} \Kres = \coprod_{\Q} F$. Moreover, by \fullref{texpre} we obtain $\Fexp=\Kexpres \preceq \Kexp$.
		
		For the converse, let $\Kexp'$ be a non-archimedean elementary extension of $\Fexp$ and let $a \in K'$ with $v(a) \neq 0$. Let $\Kexp = \dcl(F(a);\Kexp')$. By \fullref{lemexpdclres}, $\Kres = F$. By \fullref{kuhl2cor},  $v(K) \cong \coprod_{\Q} \Kres \cong G$.
	\end{proof}
	
	Now we also want to consider the exponential on a given group and give a characterisation of all countable exponential groups which can be realised as the exponential group of an ordered exponential field and a model of $\Texp$ respectively.
	
	In the following, we prove a slight enhancement of \cite[Theorem\,1.40]{kuhlmann}.
	
	\begin{lemma}\label{ctblexpgprealise}
		Let $(G,h)$ be a countable exponential group and $K$ a countable ordered field which is root closed for positive elements. Further, let $\mathbf{A}$ be a complement of $R_v$ and $\mathbf{B}$ a complement of $\U_v^{>0}$. Then $K$ can be equipped with a $v$--left exponential $\expl\co \mathbf{A} \to \mathbf{B}$ which induces the exponential group $(G,h)$ if and only if $G \cong \coprod_{\Q} \Kres$.
	\end{lemma}
	
	\begin{proof}
		Suppose that $K$ can be equipped with a $v$--left exponential which induces the exponential group $(G,h)$. Then $G$ is an exponential group in $\Kres$ and hence $G \cong \coprod_{\Q} \Kres$.
		
		Conversely, suppose that $G \cong \coprod_{\Q} \Kres$. As a consequence of the additive lexicographic decomposition properties, we have $\mathbf{A} \cong \coprod_{\Q} \Kres$ (see \cite[p.~31]{kuhlmann}). Let $\alpha\co G \to \coprod_{\Q} \Kres$ and $\beta\co \mathbf{A} \to \coprod_{\Q} \Kres$ be isomorphisms and $\tilde{\alpha}\co \Gamma \to \Q$ and $\tilde{\beta}\co G^{<0} \to \Q$ the induced isomorphisms on their underlying value sets. We define an automorphism on $\Q$ by $\tilde{\delta} = \tilde{\beta}\circ h \circ \tilde{\alpha}^{-1}$. This induces an automorphism on the ordered group $\coprod_{\Q} \Kres$ given by $\delta(s) = s\circ \tilde{\delta}^{-1}$.
		Now set $\zeta = \beta^{-1} \circ \delta \circ \alpha$. This is an isomorphism of ordered groups from $G$ to $\mathbf{A}$. The induced group exponential $\tilde{\zeta}$ is given by $\tilde{\zeta}(v_G(g)) = v(\zeta(g))$. By definition of $\tilde{\delta}$, this gives us
		$$\tilde{\zeta}(v_G(g)) = v((\beta^{-1} \circ \delta \circ \alpha)(g)) = (\tilde{\beta}^{-1} \circ \tilde{\delta} \circ \tilde{\alpha})(v_G(g)) = h(v_G(g)).$$
		
		Now $\mathbf{B}$ is isomorphic to $G$ via $-v$. Hence, $$\expl\co \mathbf{A} \to \mathbf{B}, a \mapsto ((-v)^{-1} \circ \zeta)(a) $$ defines a $v$--left exponential on $K$. It is easy to check that $\expl$ indeed induces the group exponential $h$. 
	\end{proof}
	
	We obtain strengthenings of \fullref{ansqu1} and \fullref{thmclassexpgps} for the group together with its group exponential.
	
	\begin{corollary}
		Let $G$ be a countable divisible ordered abelian group and $h$ an exponential on $G$. Then the following are equivalent.
		\begin{enumerate}
			\item $(G,h)$ can be realised as the exponential group of a non-archimedean exponential field.
			\item There exists some countable archimedean exponential field $\Fexp$ such that $G \cong \coprod_{\Q}F$.
		\end{enumerate}
	\end{corollary}
	
	\begin{corollary}\label{secondcharthm}
		Let $G$ be a countable divisible ordered abelian group and $h$ an exponential on $G$. Then the following are equivalent.
		\begin{enumerate}
			\item $(G,h)$ can be realised as the exponential group of a non-archimedean model of $\Texp$.
			\item There exists some countable $\Fexp \preceq \Rexp$ such that $G \cong \coprod_{\Q}F$, and $h$ is a strong exponential.
		\end{enumerate}
	\end{corollary}
	
	\begin{proof}
		It is a direct consequence of \fullref{thmclassexpgps} and \fullref{lemkuhl20} that (1) implies (2).
		
		Suppose that (2) holds. By \fullref{thmclassexpgps}, there exists some $\Kexp \models \Texp$ with $\Kres = F$ and $v(K) = G$. We need to equip $K$ with a $v$--left exponential such that the resulting exponential field is a model of $\Texp$ whose induced exponential group is  $(G,h)$. Note that since $\Kexp \models \Texp$, it suffices to find a (GA)--$v$--left exponential, as the $v$--middle and $v$--right exponential can stay unaltered. By \fullref{ctblexpgprealise} and \fullref{lemkuhl20}, for any decomposition of $K$, we have that $h$ induces a $v$--left exponential $\exp_h$ on $K$ satisfying (GA). Now $\tilde{\exp} = \exp_h \amalg \exp_\mathrm{M} \amalg \exp_\mathrm{R}$ gives us the required $\Texp$--exponential on $K$.
		
	\end{proof}
	
	Finally, we obtain a complete answer to \fullref{qu3} in the case that $G$ is countable.
	
	\begin{theorem}\label{thm1}
		Let $(F,G,h)$ exponential valuation triple such that $G$ is countable. Then the following are equivalent:
		\begin{enumerate}
			\item $F$ is countable, $(F,\exp) \preceq \Rexp$, $G \cong \coprod_\Q F$ and $h$ is a strong group exponential.
			\item There exists a model $\Kexp$ of $\Texp$ realising $(F,G,h)$.
		\end{enumerate}
	\end{theorem}

	\begin{remark}
		Let $\mu$ be an infinite cardinal. Then there exists a non-archimedean $\Kexp \models \Texp$ of cardinality $\mu$ such that $v(K)$ is countable if and only if $\mu \leq 2^{\aleph_0}$.
		
		Indeed, since for any valued field $K$ with value group $G$ and residue field $k$ we have $|K| \leq |k|^{|G|}$, it follows by \fullref{ctblexpgplemma2} that for any non-archimedean $\Kexp \models \Texp$ such that $v(K)$ is countable we have $|K| \leq 2^{\aleph_0}$. 
		
		Conversely, let $\mu$ be a cardinal with $\aleph_0 \leq \mu \leq 2^{\aleph_0}$ and let $\Fexp \models \Texp$ be countable and non-archimedean with value group $G$. Then it is possible to construct $\Kexp \models \Texp$ of cardinality $\mu$ which elementarily extends $\Fexp$ and has value group $G$. A construction of such a model $\Kexp$ of $\Texp$ can be deduced by applying Löwenheim--Skolem to a general construction method of uncountable o-minimal expansions of fields with countable value group which can be found in the unpublished works Fornasiero \cite{fornasiero} and Frécon \cite{frecon}.\footnote{We thank the anonymous referee for pointing out this construction method.} A more explicit construction method will be given in Krapp \cite{krapp2}.
	\end{remark}

	\section{$\kappa$--saturated Groups}\label{seckappasat}
	
	For an ordinal $\alpha$, a linearly ordered structure $(M,<,\ldots)$ is called an \textbf{$\eta_\alpha$--set} if for any $X,Y \subseteq M$ with $|X|,|Y|<\aleph_\alpha$ and $X<Y$ (i.e. $x<y$ for any $(x,y)\in X\times Y$) there exists $z \in M$ with $X<z<Y$.
	
	\textbf{Throughout this section, let $\kappa = \aleph_{\alpha}$ be an uncountable regular cardinal with $\sup_{\delta < \alpha} 2^{\aleph_\delta} \leq \aleph_\alpha$.}
	We will give answers to the questions in the introduction for the case that $G$ is a saturated group of cardinality $\kappa$.\footnote{Throughout this section, \emph{saturated} always means \emph{saturated in its own cardinality}} In particular, the results will apply if $\kappa$ is a strongly inaccessible cardinal. 
	
	It is shown in Erd\"os--Gillman--Henriksen \cite{erdos} that any two real closed fields which are $\eta_\alpha$--sets of cardinality $\aleph_\alpha$ are isomorphic. In Alling \cite{alling2} a similar result is proved for divisible ordered abelian groups, namely that any two divisible ordered abelian groups which are $\eta_\alpha$--sets of cardinality $\aleph_\alpha$ are isomorphic. Constructions of such real closed fields and divisible ordered abelian groups are given in Alling \cite{alling} and Alling--S~Kuhlmann \cite{alling3}. 
	Since both the theory of divisible ordered abelian groups and the theory of real closed fields are o-minimal, the property of being an $\eta_\alpha$--set is equivalent to $\aleph_\alpha$--saturation. Moreover, $\aleph_\alpha$--saturation of real closed fields and ordered abelian groups can be characterised in terms of their archimedean com\-po\-nents, convergence of pseudo Cauchy sequences and saturation of their value group and value set respectively (cf. F-V~Kuhlmann, S~Kuhlmann, Marshall and Zekavat \cite{kuhlmann4} and S Kuhlmann \cite{kuhlmann3}).
	
	We will exploit these properties of real closed fields and divisible ordered abelian groups which are $\eta_\alpha$--sets together with the characterisation of $\aleph_\alpha$--saturation.
	
	\begin{remark}\label{rmk:kappasat}
		\begin{enumerate}
			\item 
			It is a well-known model theoretic result that, under our given conditions on $\kappa$, for any complete theory $T$ in a countable language, there exists a saturated model of car\-di\-nal\-i\-ty $\kappa$. Moreover, a $\kappa$--saturated model is $\kappa^+$--universal. This generalises the results on real closed fields and divisible ordered abelian groups which are $\eta_\alpha$--sets as mentioned above. 
			
			\item
			Any saturated real closed field of cardinality $\kappa$ can be equipped with a $\Texp$--exponential. Indeed, if $K$ is a saturated real closed field of cardinality $\kappa$ and $\Fexp$ is a saturated model of $\Texp$ of cardinality $\kappa$, then there exists an isomorphism between $(K,+,\cdot,0,1,<)$ and $(F,+,\cdot,0,1,<)$ and $K$ can be equipped with an $\Texp$--exponential via this isomorphism.
		\end{enumerate}
	\end{remark}
	
	We will start by outlining the construction of models of $\Texp$ of cardinality $\kappa$ and value groups thereof in S~Kuhlmann--Shelah \cite{kuhlmann2}. We will also show that the constructed divisible ordered abelian groups $G_\kappa$ and the real closed fields $K_\kappa$ are saturated and of cardinality $\kappa$. Thus, by \fullref{rmk:kappasat}, all results stated for $G_\kappa$ (respectively $K_\kappa$) will apply to any saturated divisible ordered abelian group (respectively real closed field) of cardinality $\kappa$. 
	
	Note that the authors of \cite{kuhlmann2} work with logarithms rather than exponentials. Naturally, by taking inverses, in each step of the construction one can as well work with exponentials.
	Let $\Gamma_0$ be an arbitrary chain with $|\Gamma_0|\leq \kappa$. First they construct $\Gamma_\kappa$ extending $\Gamma_0$ and set $G_\kappa = (\R^{\Gamma_\kappa})_\kappa$ (cf. \cite[Section~4]{kuhlmann2}). Moreover, the construction gives us an isomorphism of ordered chains $\iota_\kappa\co \Gamma_\kappa \to G_\kappa^{<0}$, i.e. a group exponential on $G_\kappa$. 
	
	\begin{lemma} \label{lemcardkappa}
		Let $\Gamma$ be a chain with $|\Gamma| = \kappa$ and $G$ an ordered abelian group with $|G| = \kappa$. Then $|(\R^{\Gamma})_\kappa| = |\R\pow{G}_\kappa| = \kappa$.
	\end{lemma}
	
	\begin{proof}
		One easily sees that $|\R\pow{G}_\kappa|=|(\R^{\Gamma})_\kappa|$. A proof that $|(\R^{\Gamma})_\kappa| = \kappa$ is given in Alling \cite[Theorem\,2.3]{alling}. 
	\end{proof}
	
	Using \fullref{lemcardkappa} one can go through the construction of $\Gamma_\kappa$ and $G_\kappa$ to obtain $|\Gamma_\kappa| = |G_\kappa| = \kappa$. 
	
	By \cite[Proposition~4]{kuhlmann2}, there exists $\sigma \in \mathrm{Aut}(\Gamma_\kappa)$ such that $h_\kappa := \iota_\kappa \circ \sigma$ defines a strong group exponential on $\Gamma_\kappa$. Let $K_\kappa = \R\pow{G_\kappa}_\kappa$. Then $|K_\kappa| = \kappa$. By \cite[Lemma~2]{kuhlmann2}, $h_\kappa$ can be lifted to a (GA)--$v$--left exponential on $K_\kappa$. By \cite[Proposition~6]{kuhlmann2}, $K_\kappa$ can also be equipped with a restricted exponential function making it a model of $T_e$. Hence, $K_\kappa$ can be equipped with a $\Texp$--exponential $\exp$. An easy calculation shows that $(K_\kappa,\exp)$ actually induces the strong exponential group $(G_\kappa,h_\kappa)$.
	
	For an ordered group $H$ and a limit ordinal $\lambda$, a \textbf{pseudo Cauchy sequence in $H$ of legth $\lambda$} is a sequence $a_\rho \in H$, where $\rho < \lambda$, such that $v_H(a_\rho-a_\sigma)<v_H(a_\sigma-a_\tau)$ for any $\rho<\sigma<\tau<\lambda$. An element $a \in H$ is called the \textbf{pseudo limit} of $(a_\rho)_{\rho<\lambda}$ if $v_H(a_\rho-a)=v_H(a_\rho-a_{\rho+1})$ for any $\rho < \lambda$. Similarly one can define the notion of pseudo Cauchy sequences and pseudo limits for ordered fields.
	
	Using the following characterisation of saturation for divisible ordered abelian groups, we will show that $G_\kappa$ is $\kappa$--saturated (cf. S~Kuhlmann \cite[Théorème\,C]{kuhlmann3}.).
	
	\begin{theorem}
		Let $G$ be a divisible ordered abelian group and $\delta \geq 0$ an ordinal. Then $G$ is $\aleph_\delta$--saturated if and only if the following conditions hold:
		
		\begin{enumerate}\label{doagsat}
			\item $\Gamma$ is an $\eta_\delta$--set.
			\item Each archimedean component of $G$ is isomorphic to $(\R,+,0,<)$.
			\item Every pseudo Cauchy sequence in $G$ of length $\lambda < \aleph_\delta$ has a pseudo limit in $G$.
		\end{enumerate}
	\end{theorem}

	\begin{lemma}\label{etagplemma}
		Let $\Gamma$ be a $\eta_\alpha$--set of cardinality $\kappa$. Then $G=(\R^{\Gamma})_\kappa$ is a saturated divisible ordered abelian group of cardinality $\kappa$.
	\end{lemma}
	
	\begin{proof}
		By \fullref{lemcardkappa}, $|(\R^{\Gamma})_\kappa| = \kappa$. We thus only need to verify the three conditions of \fullref{doagsat}.
		
		Since $G =  (\R^{\Gamma})_\kappa$, the archimedean components of $G$ are all $\R$. Any pseudo Cauchy sequence in $G$ of length $\lambda < \kappa$ has a pseudo limit in $G$, as Hahn groups are par excellence groups in which any pseudo Cauchy sequence has a pseudo limit. In this case the restriction of the length of the pseudo Cauchy sequences ensures that the restriction on the cardinality of the support of elements in $G$ does not have an impact on this property.
	\end{proof}
	
	\textbf{In the rest of this section, we fix $\Gamma_0$ to be a chain such that $|\Gamma_0|=\kappa$ and $\Gamma_0$ is an $\eta_\alpha$--set.} The existence of such a chain is shown in Gillman \cite{gillman}. By construction, also $\Gamma_\kappa$ is a chain such that $|\Gamma_0|=\kappa$ and $\Gamma_0$ is an $\eta_\alpha$--set (and, in fact, $\Gamma_0$ and $\Gamma_\kappa$ are isomorphic as chains).
	
	\begin{proposition}\label{gkappasat}
		$G_\kappa$ is a saturated divisible ordered abelian group of cardinality $\kappa$.
	\end{proposition}
	
	\begin{proof}
		This follows directly from \fullref{etagplemma}.
	\end{proof}
	
	\begin{theorem}\label{thm2}
		Let $(F,G,h)$ be an exponential valuation triple such that $G$ is saturated and of cardinality $\kappa$. Then the following are equivalent:
		
		\begin{enumerate}
			\item $(F,G,h)$ can be realised in a model $\Kexp \models \Texp$.
			\item $F = \R$ and $h$ is a strong group exponential.
		\end{enumerate}
	\end{theorem}
	
	\begin{proof}
		Suppose that (1) holds. Since $\exp$ satisfies (GA), $h$ must be a strong group exponential. Moreover, all archimedean components of $G$ are isomorphic to $\Kres$ and to $\R$. Hence, $\Kres = \R$.
		
		For the converse, the strong group exponential $h$ gives us a strong group exponential $h_\kappa$ on $G_\kappa$ via the isomorphism $G \cong G_\kappa$. Hence, we can equip $\R\pow{G}_\kappa$ can with a $\Texp$--exponential $\exp$ such that $(\R\pow{G}_\kappa,\exp)$ realises $(\R,G,h)$. 
	\end{proof}
	
	Fix a $\Texp$--exponential $\exp$ on $K_\kappa$, whose construction was pointed out at the beginning of this section. The following theorem is due to F-V~Kuhlmann, S~Kuhlmann, Marshall and Zekavat \cite[Theorem~6.2]{kuhlmann4}.
	
	\begin{theorem}
		Let $K$ be a real closed field and $\delta \geq 0$ an ordinal. Then $K$ is $\aleph_\delta$--saturated if and only if the following conditions hold:
		
		\begin{enumerate}\label{rcfsat}
			\item $G$ is $\aleph_\delta$--saturated.
			\item $\ol{K} = \R$.
			\item Every pseudo Cauchy sequence in $K$ in a subfield of absolute transcendence degree less than $\aleph_\delta$ has a pseudo limit in $K$.
		\end{enumerate}
	\end{theorem}
	
	\begin{proposition}\label{kkappasat}
		$K_\kappa$ is $\kappa$--saturated.
	\end{proposition}
	
	\begin{proof}
		We need to check the condition from \fullref{rcfsat} for $K_\kappa = \R\pow{G_\kappa}_\kappa$. By \fullref{gkappasat}, $G_\kappa$ is $\kappa$--saturated. Clearly $\Kres = \R$. Let $F\subseteq K_\kappa$ be a subfield of $K$ of transcendence degree $\lambda < \kappa$. Then $|F| = \max\{\lambda,\aleph_0\}$. Hence, any pseudo Cauchy sequence in $F$ is either of length $\aleph_0$ or of uncountable length less than $\lambda$. Since $K_\kappa$ is a $\kappa$--bounded power series field, such pseudo Cauchy sequences have a pseudo limit in $K_\kappa$.
	\end{proof}
			
	Lastly, we will show that $G_\kappa$ has the lifting property.
	
	\begin{proposition}
		$G_\kappa$ has the lifting property.
	\end{proposition}
	
	\begin{proof}
		Since $G_\kappa = (\R^{\Gamma_\kappa})_\kappa$, the construction of a lifting of a given i\-so\-mor\-phism on $\Gamma_\kappa$ follows the same steps as the construction succeeding \fullref{kuhl2cor}.
	\end{proof}

	By \fullref{kkappasat}, we obtain the following.
	
	\begin{corollary}
		Let $G$ be the value group of an ordered exponential field $\Kexp$ such that $K$ is saturated and of cardinality $\kappa$. Then $G$ has the lifting property.
	\end{corollary}

	\section{Contraction Groups}\label{seccontr}
	
	A disadvantage of a group exponential $h$ on a group $G$ is that its domain is the value set $\Gamma$ rather than $G^{<0}$. However, composing it with $v_G$ we obtain a map $h\circ v_G\co G^{<0} \to G^{<0}$. We extend this map to a map $\chi_h$ on $G$ as follows: $$\chi_h\co G \to G, x \mapsto \begin{cases}
	(h\circ v_G)(x) &\text{ if } x<0, \\
	0 &\text{ if } x=0,\\
	(h\circ v_G)(-x) &\text{ if } x>0.
	\end{cases}$$
	$\chi_h$ is the {\bf contraction map on $G$ induced by $h$} and $(G,\chi_h)$ is the {\bf contraction group induced by $(G,h)$ 
	}.
	In general, contraction groups can be described model theoretically as follows:
	We expand the language of ordered groups to the language $\L_{\mathrm{cg}} = (+,-,0,<,\chi)$ of contraction groups, where $\chi$ is a unary function symbol. We simply write $(G,\chi)$ for an $\L_{\mathrm{cg}}$--structure. We call $(G,\chi)$ a {\bf contraction group} if it satisfies:
	\begin{itemize}
		\item $G$ is an ordered abelian group,
		\item $\chi(x) = 0$ if and only if $x = 0$,
		\item $\chi$ preserves $\leq$ (i.e. weak inequality),
		\item if $v_G(x) = v_G(y)$ and $\mathrm{sign}(x) = \mathrm{sign}(y)$, then $\chi(x) = \chi(y)$,
		\item $\chi$ is surjective.
	\end{itemize}
	It is easy to see that $(G,\chi_h)$ is indeed a contraction group.
	If $(G,\chi)$ additionally satisfies that  $x \neq 0$ implies $|x| < |\chi(x)|$, it is called \textbf{centripetal}.
	A description of contraction groups as well as some model theoretic properties such as weak o-minimality and model completeness is given in \cite[Appendix~A]{kuhlmann}.
	
	Recall that a non-archimedean exponential field $\Kexp$ induces an exponential group $(G,h_{\exp})$. We denote the corresponding contraction group by $(G,\chi_{\exp})$. Naturally, the question arises which contraction groups can be realised as the contraction groups induced by non-archimedean models of $\Texp$ (or, more generally, non-archimedean exponential fields). In this section we will use our results from the previous sections as well as \cite[Chapter~2]{kuhlmann} to give an answer for the case that $G$ is countable and the case that $G$ is $\kappa$--saturated and of cardinality $\kappa$ for some uncountable regular cardinal $\kappa = \aleph_\alpha$ with $\sup_{\delta < \alpha} 2^{\aleph_\delta} \leq \aleph_\alpha$.
	
	\begin{lemma}{\rm \cite[Lemma~2.26]{kuhlmann}} \qua\label{contrlemma1}
		Let $(G,\chi)$ be a contraction group. Let $h_\chi\co \Gamma \to G^{<0}$ be given by
		$$h_\chi(v_G(g)) = \chi(g)$$ for any $g \in G^{<0}$. Then $(G,h_\chi)$ is an exponential group. 
	\end{lemma}

	We call $h_\chi$ the {\bf group exponential induced by $\chi$}. Note that by construction $(G,h_\chi)$ induces the contraction group $(G,\chi)$.
	
	\begin{lemma}{\rm \cite[Lemma~2.28]{kuhlmann}} \qua\label{contrlemma2}
		Let $(G,h)$ be an exponential group. Then $(G,\chi_h)$ is a centripetal contraction group if and only if $(G,h)$ is a strong exponential group.
	\end{lemma}

	Using the previous two lemmas, we obtain the following versions of \fullref{thm1} and \fullref{thm2} for contraction groups.
	
	\begin{corollary}
		Let $F$ be an archimedean field, $G$ a countable divisible ordered abelian group and $\chi$ a contraction on $G$. Then the following are equivalent:
		\begin{enumerate}
			\item There exists a model $\Kexp$ of $\Texp$ such that $\Kres = F$ and $\Kexp$ induces the contraction group $(G,\chi)$.
			\item $F$ is countable, $(F,\exp) \preceq \Rexp$, $G \cong \coprod_\Q F$ and $\chi$ is centripetal.
		\end{enumerate}
	\end{corollary}

	\begin{corollary}
		Let $\kappa = \aleph_\alpha$ be an uncountable regular cardinal  with $\sup_{\delta < \alpha} 2^{\aleph_\delta} \leq \aleph_\alpha$. Let $F$ be an archimedean field, $G$ a saturated divisible ordered abelian group of cardinality $\kappa$ and $\chi$ a contraction on $G$. Then the following are equivalent:
		
		\begin{enumerate}
			\item There exists a model $\Kexp$ of $\Texp$ such that $\Kres = F$ and $\Kexp$ induces the contraction group $(G,\chi)$.
			\item $F = \R$ and $\chi$ is centripetal.
		\end{enumerate}
	\end{corollary}

	\subsection*{Acknowledgements}
		This work is part of my ongoing doctoral research project \textsl{Algebraic and Model The\-o\-ret\-ic Properties of O-minimal Exponential Fields} \cite{krapp2}. I wish to thank my doctoral supervisor Salma Kuhlmann for her continuous help and support. Many stimulating discussions with her inspired and shaped most parts of this work and her comments on previous versions of this paper helped to improve it considerably. I thank the anonymous referee for numerous helpful comments and references to known results in the literature which enabled me to rewrite and resubmit this paper.
		
		I was supported by a doctoral scholarship of Studienstiftung des deutschen Volkes as well as of Carl-Zeiss-Stiftung.

	\begin{footnotesize}
	
	\end{footnotesize}
	
\end{document}